\documentclass[reqno,11pt]{amsart}
 
\usepackage{amsmath,amsthm}
\usepackage{amssymb}
\usepackage{pb-diagram}
\usepackage{euscript}

\newcommand{\sqsp}{\renewcommand{\baselinestretch}{1.05}\tiny\normalsize}

\numberwithin{equation}{section}

\newtheorem{thm}{Theorem}[section]
\newtheorem{prop}[thm]{Proposition}
\newtheorem{cor}[thm]{Corollary}
\newtheorem{lemma}[thm]{Lemma}
\newtheorem{claim}[thm]{Claim}

\theoremstyle{definition}
\newtheorem{definition}[thm]{Definition}
\newtheorem{example}[thm]{Example}
\newtheorem{remark}[thm]{Remark}

\newcommand{\cat}[1]{{\EuScript #1}}
\newcommand{\cA}{\cat{A}}

\newcommand{\cI}{\cat{I}}

\newcommand{\bQ}{\mathbf{Q}}

\newcommand{\bZ}{\mathbf{Z}}

\newcommand{\rP}{\mathrm{P}}
\newcommand{\rSq}{\mathrm{Sq}}

\DeclareMathOperator{\Id}{Id} 
\DeclareMathOperator{\Gr}{Gr}

\begin{document}
\title[Refinement of topological realization of unstable algebras]{A $K$-theoretic refinement of topological realization of unstable algebras}
\author{Donald Yau}

\begin{abstract}
In this paper we propose and partially carry out a program to use $K$-theory to refine the topological realization problem of unstable algebras over the Steenrod algebra.  In particular, we establish a suitable form of algebraic models for $K$-theory of spaces, called $\psi^p$-algebras, which give rise to unstable algebras by taking associated graded algebras mod $p$.  The aforementioned problem is then split into (i) the \emph{algebraic} problem of realizing unstable algebras as mod $p$ associated graded of $\psi^p$-algebras and (ii) the \emph{topological} problem of realizing $\psi^p$-algebras as $K$-theory of spaces.  Regarding the algebraic problem, a theorem shows that every connected and even unstable algebra can be realized.  We tackle the topological problem by obtaining a $K$-theoretic analogue of a theorem of Kuhn and Schwartz on the so-called Realization Conjecture.
\end{abstract}


\date{\today}

\email{dyau@math.uiuc.edu}
\address{Department of Mathematics, University of Illinois at Urbana-Champaign, 1409 W. Green Street, Urbana, IL 61801 USA}

\maketitle

\sqsp

\section{Introduction}
\label{sec:intro}
Let $\cA$ be the Steenrod algebra associated to a prime $p$; the mod $p$ cohomology of a topological space is then an $\cA$-algebra.  In fact, it is an \emph{unstable} $\cA$-algebra.  This means that if $p$ is an odd prime, then $\rP^ix = x^p$ if $2i = \vert x \vert$ and $\rP^i x = 0$ if $2i > \vert x \vert$, and there are similar conditions for the prime $2$.  The problem of which unstable $\cA$-algebras (or $\cA$-modules) can be realized as the cohomology of a space has a long history and is one of the central problems in algebraic topology.  For instance, about four decades ago, Steenrod asked which polynomial algebras over $\bZ/p$ can be realized by the cohomology of a space, and the famous work of Adams and Wilkerson \cite{aw} on their embedding theorem is one of the many papers on this $\cA$-algebra realization problem.

More recently, the topological realization problem of $\cA$-modules has been studied by Kuhn and Schwartz.  In \cite{kuhn} Kuhn made an interesting conjecture, the Realization Conjecture, which says that if a finitely generated $\cA$-module is topologically realizable, then it must be finite as a set.  Thus, Kuhn's conjecture predicts that a finitely generated $\cA$-module that is not finitely generated as a $\bZ/p$-module cannot arise as the cohomology of a space.  This conjecture is partially verified by Kuhn \cite{kuhn} and is proved in its entirety by Schwartz \cite{schwartz}.  Meanwhile, Blanc \cite{blanc} approached the question of topologically realizing $\cA$-algebras from an obstruction theory point of view.  He established an obstruction theory which one can use, in principle, to decide whether any $\cA$-algebra of finite type is topologically realizable.

The purpose of this paper is to propose and partially carry out a program in which $K$-theory is used to refine the topological realization problem of unstable $\cA$-algebras.  In order to describe our program, we first need to recall a result of Atiyah relating the unstable structures of $K$-theory and mod $p$ cohomology of torsionfree spaces.

Let $X$ be a torsionfree space; that is, a space that has no torsion in integral cohomology.  Then the even dimensional part of its integral cohomology ring $H^*(X;\bZ)$ can be identified with the associated graded ring of its $K$-theory:
   \begin{equation}
   \label{eq:identification}
   \Gr^* K(X) ~\cong~ H^{\text{even}}(X;\bZ)
   \end{equation}
Here $\Gr^*(-)$ denotes the associated graded ring and the filtration on $K(X)$ arises from a skeletal filtration on $X$.  That is, letting $X_n$ denote the $n$th skeleton of $X$, the $i$th filtration ideal of $K(X)$ is the kernel $I^i(X) ~=~ \ker(K(X) \to K(X_{i-1}))$ of the restriction map.  Throughout this paper we will use $\otimes$ to denote tensor product over $\bZ$ and $\psi^n$ to denote Adams operations in $K$-theory.

Atiyah \cite[5.6 and 6.5]{atiyah} showed that, in fact, operations in $K$-theory determine those in mod $p$ cohomology.

\medskip
\begin{thm}[Atiyah]
\label{thm:atiyah}
Let $p$ be a prime and let $X$ be a torsionfree space.  If $\alpha \in K(X)$ lies in filtration $2q$, then there exist elements $\alpha_i \in K(X)$ $(i = 0, 1, \ldots, q)$ in filtration $2q + 2i(p-1)$ such that
   \begin{equation}
   \label{eq:atiyah formula}
   \psi^p(\alpha) ~=~ \sum_{i=0}^q \, p^{q-i} \alpha_i
   \end{equation}
where $\alpha_q = \alpha^p$ if $q > 0$.
This yields well-defined functions
   \[
   \label{eq:steenrod operations}
   \rP^i \colon (\Gr^* K(X)) \otimes \bZ/p \to (\Gr^{* + 2i(p-1)} K(X)) \otimes  \bZ/p,
   \]
sending $\bar{\alpha}$ $($the image of $\alpha$ in the mod $p$ associated graded ring$)$ to $\bar{\alpha}_i$.  With the identification of eq.\ \eqref{eq:identification} mod $p$, these functions $\rP^i$ are precisely the Steenrod operations $($with $\rP^i = \rSq^{2i}$ when $p = 2)$ in mod $p$ cohomology.
\end{thm}

This theorem of Atiyah suggests that there should be some suitable algebraic objects modeling $K$-theory of spaces, whose associated graded algebras mod $p$ are unstable $\cA$-algebras with Steenrod operations induced by some ``Adams operation" $\psi^p$.  As we will see shortly, there is indeed an algebraic model for $K$-theory with the desired property, and we call it $\psi^p$-\emph{algebra}.

Our program to use $K$-theory to refine the topological realization problem of unstable $\cA$-algebras can now be summarized by the following diagram.
   \begin{equation}\label{eq:program}
   \begin{diagram}
     \node[2]{(\psi^p\text{-algebras})} \arrow{s,r}{\Gr^*(-)\,\otimes\, \bZ/p} \\
     \node{\text{(spaces)}} \arrow{ne,t}{K(-)} \arrow{e,b}{H^*(-)}
       \node{(\text{unstable } \cA\text{-algebras})}
   \end{diagram}
   \end{equation}
In other words, the topological realization problem of unstable $\cA$-algebras, which corresponds to the bottom arrow, splits into two separate problems corresponding to the other two arrows:
   \begin{itemize}
   \item[(i)] The \emph{algebraic} problem of realizing unstable $\cA$-algebras as associated graded algebras mod $p$ of $\psi^p$-algebras.
   \item[(ii)] The \emph{topological} problem of realizing $\psi^p$-algebras as $K$-theory of spaces. 
   \end{itemize}
In this paper we concentrate mostly on the algebraic realization problem, but we will also present a result about the $K$-theoretic topological realization problem.

A discussion of some of the main results of this paper follows.  We begin with a result which establishes a certain kind of algebraic model for $K$-theory of spaces.  In Atiyah's Theorem \ref{thm:atiyah}, the Adams operation $\psi^p$ is the only $K$-theory operation used when relating $K$-theory and mod $p$ cohomology of torsionfree spaces.  This suggests that we should model $K$-theory by algebras having a self-map $\psi^p$ that behaves like the Adams operation on $K$-theory and that splits into a sum of the form \eqref{eq:atiyah formula} when applied to elements in the algebra.  More precisely, we define a \emph{pre}-$\psi^p$-\emph{algebra} to be a commutative, filtered ring $R$ with a distinguished endomorphism $\psi^p \colon R \to R$ such that the following two conditions hold. 
   \begin{itemize}
   \item[(i)] The $2n$th filtration of $R$ coincides with the $(2n-1)$st filtration for all $n$.
   \item[(ii)] Let $r \in R$ be in filtration $2q$ for some $q$.  If $q > 0$, then there exist elements $r_i \in R$ $(0 \leq i \leq q)$ in filtration $2q  + 2i(p-1)$ such that
   \begin{equation}\label{eq2:atiyah formula}
   \psi^p(r) ~=~ p^q r_0 + p^{q-1}r_1 + \cdots + pr_{q-1} + r_q
   \end{equation}
   where $r_q = r^p$.  If $q = 0$, then there exists an element $r^\prime \in R$ in filtration $0$ such that
   \begin{equation}\label{eq3:atiyah formula}
   \psi^p(r) ~=~ pr^\prime + r^p.
   \end{equation}
   Moreover, the elements $r_i$ are required to be well-defined in the associated graded algebra mod $p$ of $R$.
   \end{itemize}

The $K$-theory of a space together with its Adams operation $\psi^p$ is clearly a pre-$\psi^p$-algebra, at least if the space has integral cohomology that is torsionfree and is concentrated in even dimensions.  Each of equation \eqref{eq2:atiyah formula} and \eqref{eq3:atiyah formula} is called an \emph{Atiyah formula} for $r$.  Notice that the elements $r_i$ in an Atiyah formula for $r$ are not unique.

Our first main result is then the following theorem, which says that the associated graded algebra mod $p$ of a pre-$\psi^p$-algebra is ``close" to being an unstable $\cA$-algebra.

   \begin{thm}\label{thm:pre}
   Let $p$ be a fixed prime and let $R$ be a pre-$\psi^p$-algebra.  Then there exist operations
   \begin{equation}
   \label{eq:operations from R}
   \rP^i \colon \Gr^*R \otimes \bZ/p ~\to~ \Gr^{*+2i(p-1)}R \otimes \bZ/p \qquad (i \geq 0)
   \end{equation}
defined as follows.  Given any element $\overline{r} \in \Gr^{2q}R \otimes \bZ/p$ lift it to any element $r \in R$ in filtration exactly $2q$ whose image in $\Gr^{2q}R \otimes \bZ/p$ is $\overline{r}$ and write down an Atiyah formula $\psi^p(r) = \sum_{i=0}^q\, p^{q-i}r_i$ $($or, in case $q = 0$, $\psi^p(r) = r^p + pr^\prime$$)$ for $r$.  Then define $\rP^i\overline{r} \in \Gr^{2q + 2i(p-1)}R \otimes \bZ/p$ by the equation
   \begin{equation}\label{eq:def of Pi}
   \rP^i\overline{r} ~=~ \begin{cases}
       \overline{r}^p & \text{ if } i = q = 0 \\
       \overline{r}_i & \text{ if } 0 \leq i \leq q \text{ and } q > 0\\
       0              & \text{ if } i > q \geq 0. \end{cases}
   \end{equation}
   Moreover, these operations have the following four properties.
   \begin{enumerate}
   \item Each $\rP^i$ is additive.
   \item For any $q \geq 0$, $\rP^q$ is the $p$th power map from $\Gr^{2q}R \otimes \bZ/p$ to $\Gr^{2pq}R \otimes \bZ/p$.
   \item If $\vert \overline{r} \vert = 2q$ then $\rP^{i} \overline{r} = 0$ for every $i > q$.
   \item For any two elements $\overline{r}$ and $\overline{s}$ in $\Gr^*R \otimes \bZ/p$, $\rP^i(\overline{r}\overline{s}) = \sum_{l+k=i}\, (\rP^l \, \overline{r})\, (\rP^k \, \overline{s})$.
   \end{enumerate}
   \end{thm}

Thus, the associated graded algebra mod $p$ of a pre-$\psi^p$-algebra has operations $P^i$ $(i \geq 0)$ that are additive and satisfies the unstable conditions and the Cartan formula.  However, missing from the above list are the Adem relations and the basic property that $P^0 = \Id$.  Now we present some relatively simple examples of pre-$\psi^p$-algebras in which these conditions are not satisfied.

   \begin{example}\label{ex:P0}
For any prime $p$, there exists a pre-$\psi^p$-algebra structure on the ring $\bZ \lbrack \varepsilon \rbrack$ $(\varepsilon^2 = 0)$ of dual numbers, where $\varepsilon$ lies in filtration precisely $4$, in which the function $\rP^0$ as in \eqref{eq:operations from R} is not equal to the identity function.  The Adem relations, however, are satisfied for trivial reasons.
   \end{example}

   \begin{example}\label{ex:Adem}
For any prime $p > 2$, there exists a pre-$\psi^p$-algebra structure on the polynomial ring $\bZ_{(p)} \lbrack x \rbrack$ over the $p$-local integers, where $x$ lies in filtration precisely $2(p-1)$, in which the Adem relations are not satisfied.  It is the case that $P^0 = \Id$.
   \end{example}

In view of these examples, we define a $\psi^p$-\emph{algebra} to be a pre-$\psi^p$-algebra $R$ for which the operations $P^i$, $i > 0$, (as in Theorem \ref{thm:pre}) on $\Gr^*R \otimes \bZ/p$ satisfies the Adem relations and the property that $P^0 = \Id$.  We remark that it is possible to describe the Adem relations in terms of elements appearing in Atiyah formula (see Remark \ref{rem:psi p}).  A $\psi^p$-algebra map $f \colon R \to S$ between two $\psi^p$-algebras is a filtered algebra map $f \colon R \to S$ (that is, an algebra map that preserves the respective filtrations) which is compatible with the endomorphisms $\psi^p$ of $R$ and $S$, in the sense that $\psi^p f = f \psi^p$.  It follows from Atiyah's Theorem \ref{thm:atiyah} that the $K$-theory $K(X)$ with its Adams operation $\psi^p$ is a $\psi^p$-algebra, provided $X$ is a space whose integral cohomology is concentrated in even dimensions and is torsionfree.

The following result is now an immediate consequence of Theorem \ref{thm:pre}.

   \begin{cor}\label{cor:psi p}
   For any prime $p$ there is a functor 
   \[
   \Gr^*(-) \otimes \bZ/p \colon (\psi^p\text{-algebras}) \to (\text{unstable }\cA\text{-algebras})
   \]
   which associates to each $\psi^p$-algebra $R$ the unstable $\cA$-algebra $\Gr^*R \otimes \bZ/p$ with Steenrod operations $P^i$ $(i \geq 0)$ as in Theorem \ref{thm:pre}.
   \end{cor}   
Henceforth whenever $R$ is a $\psi^p$-algebra, the symbol $\Gr^*R \otimes \bZ/p$ will mean the unstable $\cA$-algebra with Steenrod operations $P^i$ $(i \geq 0)$ as in Theorem \ref{thm:pre}.

Having established that $\psi^p$-algebras are appropriate algebraic models for $K$-theory for our purpose, we now turn to the algebraic realization problem concerning $\psi^p$-algebras and unstable $\cA$-algebras.  We would like to answer the following question:  
   \begin{quote}
   Given an unstable $\cA$-algebra is there always a ``lift" of it back to a $\psi^p$-algebra?     \end{quote}
Before we can answer this question, we first need to make precise what we mean by a ``lift" of an unstable $\cA$-algebra.

For certain technical reasons, we need to work with ``connected" objects.  An unstable $\cA$-algebra $H^*$ is said to be \emph{connected} if $H^0 = \bZ/p$ and $H^{<0} = 0$.  It is called an \emph{even unstable} $\cA$-\emph{algebra} if it is trivial in odd dimensions.  Let $H^*$ be an even, connected unstable $\cA$-algebra.  A $\psi^p$-algebra $R$ is said to be \emph{connected} if the degree $0$ part $\Gr^0 R \otimes \bZ/p$ of its associated graded algebra mod $p$ is isomorphic to $\bZ/p$.  Then a \emph{lift} of $H^*$ is defined to be a connected $\psi^p$-algebra $R$ for which there exists an isomorphism
   \[
   H^* ~\cong~ \Gr^*R \otimes \bZ/p
   \]
of unstable $\cA$-algebras.

The following result shows that every even, connected unstable $\cA$-algebra has a lift that is canonical up to a choice of $\cA$-algebra generators and relations.  This result, in particular, answers the question stated above, at least if the unstable $\cA$-algebra is even and connected.

   \begin{thm}\label{thm:lift}
   Let $p$ be a fixed prime and let $H^*$ be an even, connected unstable $\cA$-algebra.  Then for each choice of a set of $\cA$-algebra generators and relations for $H^*$, there exist a corresponding canonical construction of a lift $B = B_{H^*}$ of $H^*$ and a canonical unstable $\cA$-algebra isomorphism 
   \[
   \rho_{H^*} \colon \Gr^* B_{H^*} \otimes \bZ/p ~\xrightarrow{\cong}~ H^*. 
   \]
   
Moreover, suppose that a fixed choice of $\cA$-algebra generators and relations for $H^*$ is made and that $B_{H^*}$ and $\rho_{H^*}$ are, respectively, the corresponding lift and isomorphism.  Suppose that 
   \[
   \phi \colon H^* \xrightarrow{\cong} K^*
   \]
is an isomorphism of even, connected unstable $\cA$-algebras and that the lift $B_{K^*}$ is constructed with the $\cA$-algebra generators and relations for $K^*$ corresponding via $\phi$ to those of $H^*$ already chosen.  Then there exists a canonical isomorphism 
       \begin{equation}
       \varphi \colon B_{H^*} \xrightarrow{\cong} B_{K^*}
       \end{equation}
of connected $\psi^p$-algebras such that the following diagram of unstable $\cA$-algebras and isomorphisms commutes$:$
   \begin{equation}
   \begin{diagram}
   \node{\Gr^*B_{H^*} \otimes \bZ/p} \arrow{s,l}{\rho_{H^*}} \arrow{e,t}{\varphi_*} 
      \node{\Gr^*B_{K_*} \otimes \bZ/p} \arrow{s,r}{\rho_{K^*}} \\
   \node{H^*} \arrow{e,b}{\phi} \node{K^*}
   \end{diagram}
   \end{equation}
   \end{thm}

This finishes the presentation of our main results concerning the algebraic realization problem.

As for the $K$-theoretic topological realization problem, we have a $K$-theoretic analogue of a theorem of Kuhn \cite{kuhn} and Schwartz \cite{schwartz} about the so-called Realization Conjecture that was mentioned earlier in the Introduction (see Proposition \ref{thm:K-theoretic}).  Our $K$-theoretic result is about a simpler type of algebraic objects, $\psi^p$-modules, which corresponds to $\cA$-modules.  Since we do not want to encumber this Introduction with too many definitions, we will not discuss this result in this section.

\subsection{Organization of the paper} 
This paper is organized as follows.  

In \S \ref{sec:filtered} we set the stage by discussing some preliminary materials on unstable $\cA$-algebras, filtered algebras, and pre-$\psi^p$-algebras.  Sections \ref{sec:proof of pre} contains the proofs of Theorems \ref{thm:pre}, Examples \ref{ex:P0} and \ref{ex:Adem}, and Corollary \ref{cor:psi p}.   The proof of Theorem \ref{thm:lift} is given in \S\ref{sec:proof of lift}.

The final section, \S \ref{sec:kuhn}, contains our $K$-theoretic analogue of Kuhn's Realization Conjecture.


\section{Unstable $\cA$-algebra, filtered algebra, and pre-$\psi^p$-algebra}\label{sec:filtered}

\subsection{Unstable $\cA$-algebra}\label{subsec:unstable}
Here we briefly recall the definition of an unstable algebra over the Steenrod algebra.  The reader can consult the books \cite{es}, \cite[Ch.\ 1]{schwartz2} for more information on this subject.  The field of $p$ elements is denoted by $\bZ/p$.

\subsubsection{The Steenrod algebra}\label{subsubsec:steenrod algebra}
Let $p$ be a fixed prime and denote by $\cA$ the mod $p$ Steenrod algebra.  We do not decorate $\cA$ with a subscript $p$, as some authors do, since we usually work with a fixed prime $p$ in mind.

The Steenrod algebra $\cA$ is the free, graded, associative $\bZ/p$-algebra generated by the Bockstein $\beta$ in degree $1$ and the Steenrod operations $\rP^i$ in degree $2i(p-1)$ (resp.\ $\rSq^i$ when $p = 2$) $(i \geq 0)$.  They are subject to the conditions:
   \begin{itemize}
   \item[(i)] $\rP^0 = \Id$ (resp.\ $\rSq^0 = \Id$ when $p = 2$). 
   \item[(ii)] $\beta^2 = 0$.
   \item[(iii)] The Adem relations: If $p > 2$, then for any $i, j > 0$ 
   \[
   P^iP^j ~=~ \sum_{t = 0}^{\left\lfloor i/p \right\rfloor} \, (-1)^{i+t} \binom{(p-1)(j-t) - 1}{i - pt} P^{i+j-t}P^t \qquad \text{if } i < pj, 
   \]
and
   \[
   \begin{split}
   P^i \beta  P^j &~=~ \sum_{t = 0}^{\left\lfloor i/p \right\rfloor}\, (-1)^{i+t} \binom{(p-1)(j-t)}{i-pt} \beta P^{i+j-t}P^t \\
   &\quad -\, \sum_{t=0}^{\left\lfloor (i-1)/p \right\rfloor} \, (-1)^{i+t-1} \binom{(p-1)(j-t) - 1}{i-pt-1}P^{i+j-t}\beta P^t 
   \qquad \text{if } i \leq pj.
   \end{split}
   \]
   \end{itemize}
There are similar Adem relations for the prime $2$.  (Here $\left\lfloor \frac{m}{n} \right\rfloor$ denotes the integer part of $\frac{m}{n}$, that is, the largest integer not exceeding it.)

\subsubsection{$\cA$-modules and $\cA$-algebras}\label{subsubsec:A algebra}
A module over $\cA$ is assumed to be graded by the integers.  An $\cA$-\emph{algebra} is a $\bZ/p$-algebra $M$ with an $\cA$-module structure such that the Steenrod operations satisfy the \emph{Cartan formula}.  That is, for any integer $l \geq 0$ and elements $m, m^\prime \in M$, the equality
   \[
   \rP^l(m m^\prime) ~=~ \sum_{i+j = l}\, \rP^i(m) \rP^j(m^\prime)
   \]
holds; there is a similar Cartan formula when $p = 2$.

An \emph{unstable} $\cA$-\emph{algebra} is an $\cA$-algebra $M$ which satisfies the unstable conditions: 
   \begin{itemize}
   \item[(i)] $\rP^i(m) = m^p$ if $2i = \vert m \vert$, the degree of $m$, and 
   \item[(ii)] $\beta^\varepsilon \rP^i(m) = 0$ if $2i + \varepsilon > \vert m \vert$ ($\varepsilon = 0, 1$).
   \end{itemize}
There are analogous conditions when $p = 2$.  An \emph{even unstable} $\cA$-\emph{algebra} is an unstable $\cA$-algebra $H^*$ that is concentrated in even dimensions, that is, $H^{\text{odd}} = 0$.  An unstable $\cA$-algebra is \emph{connected} if its degree $0$ part is $\bZ/p$ and $H^{<0} = 0$.

   \begin{example}
   The mod $p$ cohomology $H^*(X;\bZ/p)$ of a topological space $X$ is an unstable $\cA$-algebra.  If $X$ is a spectrum, then its mod $p$ cohomology is, in general, not an unstable $\cA$-algebra but merely an $\cA$-module.
   \end{example}


\subsection{Filtered objects}\label{subsec:filtered}

\subsubsection{Filtered group}\label{subsubsec:filtered group}
A \emph{filtered abelian group} is an abelian group $M$ together with a deceasing filtration of subgroups
   \[
   M ~=~ M^0 \supset M^1 \supset \cdots.
   \]
A map of filtered abelian groups is required to preserve the filtrations;  that is, the image of the $n$th filtration in the domain is a subset of the $n$th filtration in the target.

\subsubsection{Filtered algebra}\label{subsubsec:filtered algebra}
A \emph{filtered algebra} is a commutative ring $R$ with unit together with a multiplicative, decreasing filtration
   \[
   R = I^0 \supset I^1 \supset I^2 \supset \cdots
   \]
of ideals, called \emph{filtration ideals}.  We usually speak of a filtered algebra $R$ and leave the filtration implicit.  A map of filtered algebras is required to preserve the filtrations.

A filtration ideal $I$ of a filtered algebra $R$ is said to be \emph{closed under dividing by} $p$ if, for any element $r$ in $R$, $pr$ lies in $I$ implies that $r$ lies in $I$.  This notion will be useful when we prove the examples in the Introduction and Theorem \ref{thm:lift}.

If $R$ is a filtered algebra, then the \emph{associated graded algebra} is defined to be the (non-negatively) graded algebra
   \[
   \Gr^* R  ~\equiv~ \bigoplus_{n=0}^\infty \, I^n/I^{n+1}.
   \]
A map of filtered algebras induces naturally a map between the associated graded algebras.  If $r$ is an element in $R$, then the corresponding element in the associated graded algebra (or its tensor product with $\bZ/p$) is denoted by $\overline{r}$, and we say that $r$ is a \emph{lift} of $\overline{r}$.

   \begin{example}
One example of a filtered algebra is the complex $K$-theory of a topological space $X$ having the homotopy type of a CW complex.  The filtration ideals of $K(X)$ are the kernels
   \[
   I^n ~\equiv~ \ker(K(X) \xrightarrow{i^*} K(X_{n-1})),
   \]
where $i \colon X_{n-1} \hookrightarrow X$ is the inclusion of the $n$th skeleton.
   \end{example}

\subsection{Pre-$\psi^p$-algebra}\label{subsec:pre}

The $K$-theory $K(X)$ of a space $X$ is not just a filtered algebra; it also has $K$-theory operations.  Among these are the Adams operations $\psi^n$, $n \geq 1$, satisfying the conditions:
   \begin{itemize}
   \item[(i)] $\psi^1 = \Id$ and $\psi^n \psi^m = \psi^m \psi^n = \psi^{mn}$ for all $m$ and $n$, and
   \item[(ii)] if $p$ is a prime, $\psi^p r \equiv r^p \pmod{pK(X)}$ for all elements $r$.
   \end{itemize}

As discussed in \S \ref{sec:intro}, Atiyah's Theorem \ref{thm:atiyah} says that for a space whose integral cohomology is torsionfree and is concentrated in even dimensions, the Adams operation $\psi^p$ on its $K$-theory determines the Steenrod operations $P^i$ $(i \geq 0)$ on its mod $p$ cohomology.  Since we would like to model $K$-theory of spaces by algebraic objects which corresponds to unstable $\cA$-algebras, as $K$-theory corresponds to mod $p$ cohomology, Atiyah's result suggests that we consider the following.

We already stated the definition of a pre-$\psi^p$-algebra in \S \ref{sec:intro}, but we will repeat it here.   
   \begin{definition}\label{def:pre}
   Let $p$ be an arbitrary but fixed prime.  Define a \emph{pre}-$\psi^p$-\emph{algebra} to be a filtered algebra $R$ that comes equipped with a distinguished endomorphism $\psi^p \colon R \to R$ such that the following two conditions hold: 
   \begin{itemize}
   \item[(i)] The $2n$th filtration of $R$ coincides with the $(2n-1)$st filtration for all $n$. 
   \item[(ii)] Let $r \in R$ be in filtration $2q$ for some $q$.  If $q > 0$, then there exist elements $r_i \in R$ $(0 \leq i \leq q)$ in filtration $2q  + 2i(p-1)$ such that
   \begin{equation}\label{eq2':atiyah formula}
   \psi^p(r) ~=~ p^q r_0 + p^{q-1}r_1 + \cdots + pr_{q-1} + r_q
   \end{equation}
   where $r_q = r^p$.  If $q = 0$, then there exists an element $r^\prime \in R$ in filtration $0$ such that
   \begin{equation}\label{eq3':atiyah formula}
   \psi^p(r) ~=~ pr^\prime + r^p.
   \end{equation}
   Moreover, the elements $r_i$ in \eqref{eq2':atiyah formula} are required to be well-defined in the associated graded algebra mod $p$ of $R$.
   \end{itemize}
   \end{definition}
Equation \eqref{eq2':atiyah formula} (or \eqref{eq3':atiyah formula}, if $q = 0$) is called an \emph{Atiyah formula} for $r$.  Note that the elements $r_i$ in an Atiyah formula for $r$ are not unique.

   \begin{remark}\label{rem1:pre}
   We now pause to make a few remarks about this definition. 
   \begin{itemize}
   \item[(i)] An element $r \in R$ can be considered to lie in different filtrations, since if $r$ lies in filtration $2q \geq 2$ then it also lies in filtration $2(q - 1)$.  Thus, condition (ii) above has to be interpreted to mean that for every element $r \in R$, regardless of what filtration (say, $2q$) it is considered to be in, $\psi^p(r)$ can be decomposed into the form \eqref{eq2':atiyah formula} or \eqref{eq3':atiyah formula} with the $r_i$ well-defined in the associated graded algebra mod $p$.  
   \item[(ii)] If $r$ lies in filtration $2q$ with $q > 0$, then an Atiyah formula for $r$ also yields an Atiyah formula when $r$ is considered to be in filtration $2(q-1)$, since we can rewrite \eqref{eq2':atiyah formula} as
   \[
   \psi^p(r) ~=~ p^{q-1}(pr_0) + \cdots + p^2(pr_{q-3}) + p(pr_{q-2} + r_{q-1}) + r^p.
   \]
   \item[(iii)] If we think of an element $r \in R$ to lie in filtration, say, $2q$, then $\psi^p(r)$ is considered to lie in the same filtration.  In particular, if $q > 0$ then 
   \[
   \overline{\psi^p(r)} ~=~ 0
   \]
in $\Gr^{2q} R \otimes \bZ/p$.
   \item[(iv)] The requirement that the $r_i$ be well-defined in the associated graded algebra mod $p$ means that if $\psi^p(r)$ admits another decomposition of the form \eqref{eq2':atiyah formula}, say, 
   \[
   \psi^p(r) ~=~ \sum_{i=0}^q p^{q-i} r_i^\prime, 
   \]
then 
   \[
   \overline{r_i} ~=~ \overline{r_i}^\prime
   \]
in $\Gr^{2q + 2i(p-1)}R \otimes \bZ/p$ for all $i$. 
   \end{itemize}
   \end{remark}

The $K$-theory of a space together with its Adams operation $\psi^p$ is clearly a pre-$\psi^p$-algebra, at least if the space has integral cohomology that is torsionfree and is concentrated in even dimensions.

If $R$ and $S$ are pre-$\psi^p$-algebras, then a \emph{pre}-$\psi^p$-\emph{algebra map} 
   \[
   f \colon R \to S
   \]
is a map of the underlying filtered algebras that is compatible with $\psi^p$, in the sense that 
   \[
   f \psi^p ~=~ \psi^p f.  
   \]
An ideal $I \subset R$ is said to be a $\psi^p$-\emph{ideal} if it is closed under $\psi^p$.  If $I \subset R$ is a $\psi^p$-ideal and that $R/I$ inherits from $R$ the structure of a pre-$\psi^p$-algebra, then the natural map 
   \[
   R \twoheadrightarrow R/I
   \]
is a pre-$\psi^p$-algebra map.

Now we make a few preliminary observations concerning Atiyah formula; they will be used below to show that certain filtered algebras with a distinguished endomorphism $\psi^p$ admit pre-$\psi^p$-algebra structures.  The first preliminary result says that if two elements ``admit Atiyah formulas," then so does their sum.

   \begin{lemma}
   \label{lem:Atiyah}
    Let $p$ be a prime and let $R = (R, \lbrace I^n \rbrace)$ be a filtered algebra with a distinguished endomorphism $\psi^p$ and in which the $2n$th filtration coincides with the $(2n-1)$st filtration for all $n$.  Suppose that $r$ and $s$ are elements in $R$ with $r \in I^{2q} \setminus I^{2q+2}$ and $s \in I^{2v} \setminus I^{2v+2}$ for some integers $q < v$.  Then:
   \begin{itemize}
   \item[(i)] If $q = 0$ and if $\psi^p(r)$ and $\psi^p(s)$ admit decompositions of the forms, respectively, \eqref{eq3':atiyah formula} and \eqref{eq2':atiyah formula} $($where the $\overline{s_i}$ are not-necessarily well-defined$)$ then $\psi^p(r + s)$ also admits a decomposition of the form \eqref{eq3':atiyah formula}.
   \item[(ii)] If $q > 0$ and if both $\psi^p(r)$ and $\psi^p(s)$ admit decompositions of the form \eqref{eq2':atiyah formula} $($where the $\overline{r_i}$ and $\overline{s_i}$ are not-necessarily well-defined$)$, then so does $\psi^p(r + s)$.
   \end{itemize}
   \end{lemma}

   \begin{proof}
The proofs for the two statements are similar, with the first one being easier, so we will only present the argument for the second statement.

Let us write $t = r + s$ and note that $t$ lies in $I^{2q} \setminus I^{2q+2}$.  By assumption we can write  
   \[
   \psi^p(r) ~=~ \sum_{i=0}^q p^{q-i} r_i 
   \]
and
   \[
   \psi^p(s) ~=~ \sum_{i=0}^v p^{v-i} s_i.
   \]
Define the following elements
   \[
   t_i ~=~ \begin{cases}
                r_0 + s^\prime        & \text{ if } i = 0             \\
                r_i + s_i             & \text{ if } 1 \leq i \leq q-2 \\
                r_{q-1} + s_{q-1} - c & \text{ if } i = q-1           \\
                (r + s)^q             & \text{ if } i = q.
                  \end{cases}
   \]
Here $s^\prime$ and $c$ are given by
   \[
   \begin{split}
   s^\prime & ~=~ \sum_{i=0}^{v-q}\, p^{v-q-i}s_i  \\
   c        & ~=~ \sum_{i=1}^{p-1}\, \frac{1}{p} \binom{p}{i}r^{p-i}s^i. 
   \end{split}
   \]
Note that $c$ satisfies the equation
   \[
   (r + s)^p - pc ~=~ r^p + s^p.
   \]
Then we have that
   \begin{equation*}
   \label{eq:r+s}
   \begin{split}
   \psi^p(t) & ~=~  \sum_{i=0}^q p^{q-i} r_i ~+~ \sum_{i=0}^v p^{v-i} s_i \\
             & ~=~ p^q(r_0 + s^\prime) ~+~ \sum_{i = 1}^{q-2}\, p^{q-i} (r_i + s_i) \\
             & \qquad \qquad ~+~ p(r_{q-1} ~+~ s_{q-1} - c) ~+~ (r + s)^p \\
             & ~=~ \sum_{i=0}^q p^{q-i}t_i.
   \end{split}
   \end{equation*}

The Lemma now follows.
\end{proof}

The following lemma gives a sufficient condition in order that a filtered algebra with a distinguished endomorphism $\psi^p$ be a pre-$\psi^p$-algebra.  Recall that a filtration ideal $I$ in a filtered algebra $R$ is \emph{closed under dividing by} $p$ if, for any element $r$ in $R$, $pr \in I$ implies that $r \in I$.

\begin{lemma}\label{lem:dividing by p}
Let $p$ be a prime and let $R$ be a filtered algebra with a distinguished endomorphism $\psi^p$ in which every filtration ideal is closed under dividing by $p$.  Assume that for some element $r \in R$ in filtration $2q$ with $q > 0$, $\psi^p(r)$ admits two decompositions of the form \eqref{eq2':atiyah formula}$:$
   \begin{equation}\label{eq:assumption}
   \psi^p(r) ~=~ \sum_{i=0}^q\, p^{q-i} r_i ~=~ \sum_{i=0}^q\, p^{q-i} r_i^\prime.
   \end{equation}
Then 
   \[
   \overline{r_i} ~=~ \overline{r_i}^\prime
   \]
in $\Gr^{2q + 2i(p-1)}R \otimes \bZ/p$ for $0 \leq i \leq q$.
\end{lemma}

\begin{proof}
Let $I^n$ denote the $n$th filtration ideal of $R$.  It follows from \eqref{eq:assumption} that we have, modulo filtration $2q + 2(p-1)$, the equality 
   \[
   p^q(r_0^\prime - r_0) ~=~ 0.
   \]
Since $I^{2q+2(p-1)}$ is closed under dividing by $p$, the above equality implies that
   \[
   r_0^\prime - r_0 ~=~ s_0
   \]
for some element $s_0 \in I^{2q+2(p-1)}$.  In particular, the Lemma is true in the case $i = 0$.

We will now show by induction that for $1 \leq k \leq q-1$, there exist elements $s_k$ in filtration $2q + 2(k+1)(p-1)$ such that
   \begin{equation}\label{eq:r_k}
   r_k^\prime - r_k ~=~ s_k - ps_{k-1},
   \end{equation}
in which $s_0$ is as in the previous paragraph.  For the case $k = 1$, \eqref{eq:assumption} implies that, modulo filtration $2q + 2(2)(p-1)$, the equality
   \[
   p^q r_0 + p^{q-1}r_1 ~=~ p^q r_0^\prime + p^{q-1}r_1^\prime
   \]
holds, and so by rearranging terms we see that
   \[
   0 ~=~ p^{q-1} (ps_0 + r_1^\prime - r_1).
   \]
Since $I^{2q+2(2)(p-1)}$ is closed under dividing by $p$, it follows that
   \[
   ps_0 + r_1^\prime - r_1 ~=~ s_1
   \]
for some element $s_1 \in I^{2q+2(2)(p-1)}$.

Suppose that \eqref{eq:r_k} has been proved for all $k$, $1 \leq k \leq n$, for some $n < q-1$.  Then, working modulo filtration $2q + 2(n+2)(p-1)$, \eqref{eq:assumption} implies that
   \[
   \begin{split}
   0 &~=~ \sum_{i=0}^{n+1} p^{q-i} (r_i^\prime - r_i) \\
     &~=~ \sum_{i=0}^{n} p^{q-i} (s_i - ps_{i-1}) ~+~ p^{q-n-1}(r_{n+1}^\prime - r_{n+1}) \qquad (s_{-1} \equiv 0) \\
     &~=~ p^{q-n}s_n +  p^{q-n-1}(r_{n+1}^\prime - r_{n+1}) \\
     &~=~ p^{q-n-1}(ps_n + r_{n+1}^\prime - r_{n+1}). 
   \end{split}
   \]
Thus, since $I^{2q + 2(n+2)(p-1)}$ is closed under dividing by $p$, the case $k = n+1$ follows.  Therefore, \eqref{eq:r_k} holds for all $k$, $1 \leq k \leq q -1$.

In particular, since the case $i = q$ is trivially true (because $r_q = r^p = r_q^\prime$), we see that
   \[
   \overline{r_i} ~=~ \overline{r_i}^\prime
   \]
for $0 \leq i \leq q$, thereby proving the Lemma. 
\end{proof}


\section{Proof of Theorem \ref{thm:pre}, Examples \ref{ex:P0} and \ref{ex:Adem}, and Corollary \ref{cor:psi p}}\label{sec:proof of pre}

We begin with the proof of Theorem \ref{thm:pre}.

\subsection{Proof of Theorem \ref{thm:pre}}\label{subsec:proof of thm:pre}
First we prove that the operations $P^i$, as defined in eq.\ \eqref{eq:def of Pi}, are well-defined.  We will use the notation in Theorem \ref{thm:pre}, so $\overline{r}$ is an element in $\Gr^{2q} R \otimes \bZ/p$.  The cases $q = 0$ and $q > 0$ need to be considered separately.  However, the arguments are very similar, so we will only present the argument for the more complicated case $q > 0$.

For this case, it suffices to prove the following statement:  
   \begin{quote}
Let $s$ be another lift of $\overline{r}$ back to $R$, so that 
   \[
   s ~=~ r + ph + f
   \]
   for some elements $h$ and $f$ in filtrations, respectively, $2q$ and $2q + 2n$ for some $n \geq 1$.  Then there exists an Atiyah formula 
   \[
   \psi^p(s) ~=~ p^qs_0 + p^{q-1}s_1 + \cdots + ps_{q-1} + s_q
   \]
for $s$ (in filtration $2q$) such that
   \[
   \overline{r}_i ~=~ \overline{s}_i
   \]
in $\Gr^{2q + 2i(p-1)}R \otimes \bZ/p$ for $i = 0, \ldots , q$.
   \end{quote}

To begin proving this statement, let us first write down Atiyah formulas for the elements $h$ and $f$,
   \[
   \begin{split}
   & \psi^p(h) ~=~ p^qh_0 + p^{q-1}h_1 + \cdots + ph_{q-1} + h^p \\
   & \psi^p(f) ~=~ p^{q+n}f_0 + p^{q+n-1}f_1 + \cdots + pf_{q+n-1} + f^p,
   \end{split}
   \]
with $h_i$ in filtration $2(q + i(p-1))$ and $f_i$ in filtration $2(q + n + i(p-1))$.  Define an element $\gamma$ in $R$ by the equation
   \[
   r^p + ph^p + f^p ~=~ s^p + p\gamma ~=~ (r + ph + f)^p + p\gamma,
   \]
and define elements $s_i$ as follows:
   \[
   s_i ~=~   
     \begin{cases}
     r_i + ph_i + p^n f_i & \text{ if} \quad 0 \leq i \leq q-2 \\
     r_{q-1} + ph_{q-1} + \gamma + \sum_{j=q-1}^{q+n-1} p^{q+n-j-1} f_j  & \text{ if} \quad  i = q-1 \\
     s^p & \text{ if} \quad i = q.
     \end{cases}
   \]
Now we calculate
   \[
   \begin{split}
   \psi^p(s)
   & ~=~ \psi^p(r) + p\psi^p(h) + \psi^p(f) \\
   & ~=~ \biggl(\sum_{i = 0}^{q-1}p^{q-i} r_i\biggr) + p \biggl(\sum_{i=0}^{q-1}p^{q-i}h_i\biggr) \\
   & \qquad + \biggl(\sum_{i=0}^{q+n-1} p^{q+n-i} f_i\biggr) + r^p + ph^p + f^p \\
   & ~=~ \biggl(\sum_{i = 0}^{q-1} p^{q-i} r_i\biggr) + p \biggl(\sum_{i=0}^{q-1}p^{q-i}h_i\biggr) + \biggl(\sum_{i=0}^{q+n-1} p^{q+n-i} f_i\biggr) + s^p + p\gamma \\
   & ~=~ \sum_{i=0}^{q-2} p^{q-i}\left(r_i + ph_i + p^n f_i\right) \\
   & \qquad + p\biggl(r_{q-1} + ph_{q-1} + \gamma + \sum_{j=q-1}^{q+n-1} p^{q+n-j-1} f_j \biggr) + s^p \\
   & ~=~ \sum_{i=0}^{q}p^{q-i}s_i.
   \end{split}
   \]
It is not hard to see that the elements $s_i$ satisfy the required properties.  For instance,    
   \[
   \overline{s}_{q-1} ~=~ \overline{r}_{q-1} 
   \]
because $\gamma$ lies in filtration $2pq$, $f_j$ (for $j \geq q-1$) lies in filtration $2(q + n + (q-1)(p-1))$, and $ph_{q-1}$ is $p$-divisible.

Thus, the functions $P^i$ as defined in eq.\ \eqref{eq:def of Pi} are well-defined.

Now we consider the four properties of the $P^i$ listed in Theorem \ref{thm:pre}.  The first three of them are immediate consequences of the definition (eq.\ \eqref{eq:def of Pi}) of the functions $P^i$.

As for the last statement, let $\overline{r}$ and $\overline{s}$ be elements of $\Gr^*R \otimes \bZ/p$ in degrees $2m$ and $2n$, respectively.  Without loss of generality we may assume that $m \leq n$.  We will denote by $r$ and $s$ arbitrary but fixed lifts of $\overline{r}$ and $\overline{s}$, respectively, to $R$ in filtrations precisely $2m$ and $2n$.  Note that $rs$ is a lift of $\overline{r}\overline{s}$.

The case when both $m$ and $n$ are equal to $0$ is immediate.  So let us now consider the case when $m = 0$ and $n > 0$.  We write down Atiyah formulas for $r$ and $s$: 
   \[
   \begin{split}
   \psi^p(r) & ~=~ r^p + pr^\prime \\
   \psi^p(s) & ~=~ p^n s_0 ~+~ \cdots ~+~ ps_{n-1} + s^p. 
   \end{split}
   \]
Here $r^\prime$ is some element in $R$ in filtration $0$.  Using the fact that $\psi^p$ is multiplicative, we have that
   \[
   \begin{split}
   \psi^p(rs) & ~=~ (r^p + pr^\prime)\, \sum_{i=0}^n\, p^{n-i}s_i  \qquad (s_n \equiv s^p) \\
              & ~=~ \biggl( \sum_{i=0}^{n-1}\, p^{n-i} s_i (r^p + pr^\prime) \biggr) ~+~ s^p(r^p + pr^\prime) \\
   \end{split}
   \]
Since 
   \[
   \rP^i\overline{r} ~=~ \begin{cases}
      \overline{r}^p & \text{ if } i = 0 \\
      0              & \text{ if } i > 0, \end{cases}
   \]
the last statement of the theorem when $m = 0$ and $n > 0$ follows.

Finally\label{page:cartan}, we consider the case when both $m$ and $n$ are positive.  The Atiyah formula for $s$ is as above, but that for $r$ looks like
   \[
   \psi^p(r) ~=~ p^m r_0 ~+~ \cdots ~+~ p r_{m-1} ~+~ r^p.
   \]
Therefore, we have that 
   \[
   \psi^p(rs) ~=~ \sum_{i=0}^{m+n}\, p^{m+n-i}c_i,
   \]
where
   \[
   c_i ~=~ \sum_{l+k = i}\, r_l s_k.
   \]
The case when $m, n > 0$ follows.  This finishes the proof of statement (iv) in the theorem.

The proof of Theorem \ref{thm:pre} is complete.

\begin{remark}\label{rem:product}
The proof of statement (iv) in Theorem \ref{thm:pre} above also shows that if two elements ``admit Atiyah formulas," then so does their product.  This is the multiplicative analogue of Lemma \ref{lem:Atiyah}.
\end{remark}

\subsection{Proof of Examples \ref{ex:P0} and \ref{ex:Adem}}\label{subsec:proof of ex:P0}

Now we prove the two examples in the Introduction.

\begin{proof}[Proof of Example \ref{ex:P0}]
Fix a prime $p$ and let $R$ be the filtered ring $\bZ \lbrack \varepsilon \rbrack$ $(\varepsilon^2 = 0)$ of dual numbers with the $\varepsilon$-adic filtration, where $\varepsilon$ lies in filtration precisely $4$.  Let $k$ be any integer and define the endomorphism $\psi^p_k$ on $R$ by
   \begin{equation}\label{eq1:P0}
   \psi^p_k(\varepsilon) ~=~ p^2 k\varepsilon.
   \end{equation}
Using Lemmas \ref{lem:Atiyah} and \ref{lem:dividing by p} it is readily checked that $R$ is a pre-$\psi^p$-algebra (Definition \ref{def:pre}) with $\varepsilon_0 = k\varepsilon$.  Thus, in $\Gr^4 R \otimes \bZ/p$ we have    
   \[
   \rP^0\overline{\varepsilon} ~=~ k\overline{\varepsilon},
   \]
which is equal to $\overline{\varepsilon}$ if and only if 
   \[
   k \equiv 1 \pmod p.
   \]
In other words, this congruence condition is equivalent to the condition that $\rP^0 = \Id$.  In particular, if $k \not\equiv 1 \pmod p$, then $\rP^0 \not= \Id$.

It is worth pointing out that the Adem relation \emph{is} satisfied in $\Gr^*R \otimes \bZ/p$, since only $\rP^0$ can be non-zero.
\end{proof}

                 
\begin{proof}[Proof of Example \ref{ex:Adem}]\label{subsec:proof of ex:Adem}
Fix a prime $p > 2$ and let $\bZ_{(p)}$ denote the ring of integers localized at $p$.  Let $R$ be the filtered polynomial ring $\bZ_{(p)} \lbrack x \rbrack$ with the $x$-adic filtration, where $x$ lies in filtration precisely $2(p-1)$.  Define the endomorphism $\psi^p$ on $R$ by the equation
   \[ 
   \begin{split}
   \psi^p(x) & ~=~ -p^{p-2}x^2 + \sum_{i=1}^p p^{p-i} x^i \\
             & ~=~ p^{p-1}x + p^{p-3}x^3 + \cdots + px^{p-1} + x^p.
   \end{split}
   \]
It follows from Lemmas \ref{lem:Atiyah} and \ref{lem:dividing by p} and Remark \ref{rem:product} that $R$ is a pre-$\psi^p$-algebra.

Now the operation $\rP^i$ $(0 \leq i \leq p - 1)$ takes $\overline{x} \in \Gr^{2(p-1)}R \otimes \bZ/p$ to
   \[
   \rP^i\overline{x} ~=~ \begin{cases} 
       \overline{x}^{i+1} & \text{if } i \not= 1 \\
       0                  & \text{if } i = 1. \end{cases}
   \]
In particular, we have that
   \[
   \rP^1 \rP^1 \overline{x}  ~=~ \rP^1 0  ~=~ 0,
   \]
which is \textbf{not} equal to 
   \[
   2\rP^2 \overline{x} ~=~ 2\overline{x}^3,
   \]
since $p > 2$.  In other words, $\rP^1 \rP^1 \not= 2\rP^2$.

In summary, $R$ is a pre-$\psi^p$-algebra for which the operations $\rP^i$ on $\Gr^* R \otimes \bZ/p $ do not satisfy the Adem relation $\rP^1 \rP^1 = 2\rP^2$.

This finishes the proof of Example \ref{ex:Adem}.
\end{proof}


\subsection{Proof of Corollary \ref{cor:psi p}}\label{subsec:proof of cor psi p}

In view of Examples \ref{ex:P0} and \ref{ex:Adem}, to make the associated graded algebra mod $p$ of a pre-$\psi^p$-algebra into an unstable $\cA$-algebra, we should add the assumptions that the operations $P^i$ satisfy $P^0 = \Id$ and the Adem relations.  Thus, we make the following definition.

   \begin{definition}\label{def:psi p}
   A $\psi^p$-\emph{algebra} is a pre-$\psi^p$-algebra $R$ for which the operations $P^i$ on the associated graded algebra mod $p$ (as in Theorem \ref{thm:pre}) satisfy the Adem relations and the property $P^0 = \Id$.  A $\psi^p$-\emph{algebra map} is a map of the underlying pre-$\psi^p$-algebras.  A $\psi^p$-algebra $R$ is \emph{connected} if there is an isomorphism    
   \[
   \Gr^0 R \otimes \bZ/p ~\cong~ \bZ/p.
   \]
   \end{definition}

   \begin{remark}\label{rem:psi p}
   We pause to make a few comments about this definition.
   \begin{itemize}
   \item[(i)] It follows from Theorem \ref{thm:pre} that the associated graded algebra mod $p$, $\Gr^*R \otimes \bZ/p$, of a $\psi^p$-algebra $R$ is an unstable $\cA$-algebra in which the Steenrod operations are given by \eqref{eq:def of Pi}.
   \item[(ii)] Atiyah's Theorem \ref{thm:atiyah} implies that the $K$-theory $K(X)$ with its Adams operation $\psi^p$ is a $\psi^p$-algebra, provided that $X$ is a space whose integral cohomology is concentrated in even dimensions and is torsionfree.  
   \item[(iii)] Examples \ref{ex:P0} and \ref{ex:Adem} show that there are  pre-$\psi^p$-algebras that are not $\psi^p$-algebras.
   \item[(iv)] The condition that the operations $P^i$ satisfy the Adem relations can be rephrased in terms of elements appearing in Atiyah formulas, as follows: For each element $r \in R$ in filtration, say, $2q$ with $q > 0$, there exist elements $r_i$ ($0 \leq i \leq q$) and $r_{i,j}$ $(0 \leq j \leq q + i(p-1))$ in $R$ in filtrations, respectively, $2q + 2i(p-1)$ and $2q + 2(i+j)(p-1)$ such that
     \begin{equation}\label{eq:r_i,j}
       \begin{split}
     \psi^p(r)    &~=~ \sum_{i=0}^q\, p^{q-i} r_i \\
     \psi^p(r_i)  &~=~ \sum_{j=0}^{q+i(p-1)} p^{q+i(p-1)-j} r_{i,j}
       \end{split}
     \end{equation}
These elements are required to satisfy the following condition.   Whenever $i, j > 0$ and $i < pj$, the following equality holds in $\Gr^{2(q + (i+j)(p-1))}R \otimes \bZ/p:$
     \begin{equation}\label{eq:Adem}
       \begin{split}
       \overline{r}_{j,i} & ~=~
       \sum_{t=0}^{\left\lfloor i/p \right\rfloor}\, (-1)^{i+t} \binom{(p-1)(j-t)-1}{i-pt} \overline{r}_{t,i+j-t} \quad \text{ if} \quad  p > 2 \\
       \overline{r}_{j,i} & ~=~
       \sum_{t=0}^{\left\lfloor i/2 \right\rfloor}\, \binom{2j - 2t - 1}{2i - 4t} \overline{r}_{t,i+j-t} \quad \text{ if} \quad p = 2. 
       \end{split}
     \end{equation}
   \item[(v)] If 
      \[
      f \colon R \to S
      \]
   is a pre-$\psi^p$-algebra map, then the induced map on the associated graded algebras mod $p$,
     \[
     f_* \colon \Gr^*R \otimes \bZ/p \to \Gr^*S \otimes \bZ/p,
     \]
commutes with the operations $P^i$.  In particular, this applies when $f$ is a map of $\psi^p$-algebras.
   \end{itemize}
   \end{remark}

\begin{proof}[Proof of Corollary \ref{cor:psi p}]
This is an immediate consequence of Theorem \ref{thm:pre}, Definition \ref{def:psi p}, and Remark \ref{rem:psi p} (v).
\end{proof}


\section{Proof of Theorem \ref{thm:lift}}\label{sec:proof of lift}

Before giving the proof, we first state the definition of a lift of an even and connected unstable $\cA$-algebra.

   \begin{definition}\label{def:lift}
   Let $p$ be a prime and let $H^*$ be an even and connected unstable $\cA$-algebra.  Define a \emph{lift} of $H^*$ to be a connected $\psi^p$-algebra (see Definition \ref{def:psi p}) $R$ for which there exists an isomorphism
   \[
   H^* ~\cong~ \Gr^*R \otimes \bZ/p
   \]
of unstable $\cA$-algebras.  
   \end{definition}

We are now ready for the proof of Theorem \ref{thm:lift}.

\subsection{Construction of $B_{H^*}$}

\subsubsection{Choosing generators and relations for $H^*$}
We begin the proof with the construction of the lift $B_{H^*}$ of $H^*$.

Let us choose a set of $\cA$-algebra generators $\lbrace x_\theta \rbrace_{\theta \in \Theta}$ and relations $\lbrace f_\delta \rbrace_{\delta \in \Delta}$ (all in positive, even dimensions) for the connected, even unstable $\cA$-algebra $H^*$.  Denote the degree $\vert x_\theta \vert$ of $x_\theta$ by $2d_\theta$ for some positive integer $d_\theta$.  Thus, every element in $H^*$ can be written as a polynomial, with coefficients in $\bZ/p$, in elements of the form
   \begin{equation}\label{eq:iterated Pi}
   \rP^{\mathbf{i}}x_\theta ~=~ \rP^{i_n} \cdots \rP^{i_1}x_\theta.
   \end{equation}
Here $n \geq 0$, $\theta \in \Theta$, and, for $0 \leq k \leq n-1$, $i_{k+1}$ is in the range
   \begin{equation}\label{eq:range of i}
   0 \leq i_{k+1} \leq d_\theta + (p-1) \sum_{j=1}^k\, i_j.
   \end{equation}
The relations among the $x_\theta$ and the Adem relations (when applied to these elements) are contained in the set $\lbrace f_\delta \rbrace$.  Denoting the ideal in $H^*$ generated by the $f_\delta$ by $\langle f_\delta \colon \delta \in \Delta \rangle$, we have the identification
   \begin{equation}\label{eq:rep of H}
   H^* ~=~ \frac{\bZ/p \left\lbrack \, \lbrace \rP^{i_n} \cdots \rP^{i_1}x_\theta \,
   \colon \theta \in \Theta, \, n \geq 0 \rbrace \, \right\rbrack}{\langle f_\delta \colon \delta \in \Delta \rangle}
   \end{equation}
in which the $i_{k+1}$ are in the range \eqref{eq:range of i}.

\subsubsection{Construction of $B_{H^*}$}
The desired lift $B_{H^*}$ is going to be a quotient of the filtered polynomial algebra
   \[
   \Pi ~\buildrel \text{def} \over =~ \bZ \left \lbrack \, \lbrace X_{(\theta, i_1, \ldots , i_n)} \rbrace \, \right\rbrack.
   \]
Here $n \geq 0$, $\theta \in \Theta$, and the $i_k$ are in the range \eqref{eq:range of i}.  The indeterminate $X_{(\theta, i_1, \ldots , i_n)}$ is in filtration exactly 
   \[
   2d_\theta + 2(p-1)\sum_{j=1}^n i_j,
   \]
which is equal to the degree of the element in \eqref{eq:iterated Pi}.  Define an endomorphism $\psi^p$ on $\Pi$ as follows:  Given a multi-index $(\theta, \mathbf{i}) = (\theta, i_i, \ldots, i_n)$, let 
   \[
   \sigma ~=~ d_\theta + (p-1) \sum_{j=1}^n\, i_j
   \]                                                                                  
and define
   \begin{equation}\label{eq:psi p X}
   \psi^p X_{(\theta, \mathbf{i})} ~=~ p^\sigma X_{(\theta, \mathbf{i})} +
    \biggl( \sum_{i=1}^{\sigma - 1} \, p^{\sigma - i} X_{(\theta, \mathbf{i}, i)} \biggr) + X_{(\theta, \mathbf{i})}^p.
   \end{equation}
As in Examples \ref{ex:P0} and \ref{ex:Adem}, using Lemmas \ref{lem:Atiyah} and \ref{lem:dividing by p} and Remark \ref{rem:product}, it is readily checked that $\Pi$ together with the endomorphism $\psi^p$ is a pre-$\psi^p$-algebra.

Now for each $\delta \in \Delta$ let $f_\delta$ also denote the element in $\Pi$ obtained canonically from the $f_\delta$ in $H^*$ by replacing each nonzero coefficient $a \pmod p$ by the unique positive integer $a$ not exceeding $p-1$ and the element in $\eqref{eq:iterated Pi}$ by $X_{(\theta, i_1, \ldots, i_n)}$.  The $k$-fold composite, for any non-negative integer $k$, of $\psi^p$ is denoted $\psi^{pk}$, where $\psi^0$ is defined to be the identity map.  Let $\cI \subset \Pi$ be the ideal generated by the set of elements
   \[
   \left\lbrace \psi^{pk}f_\delta \, \colon \, k \geq 0, \, \delta \in \Delta \right\rbrace.
   \]
Then the proposed lift $B_{H^*}$ is defined to be the quotient
   \begin{equation}\label{eq:def of B}
   B ~=~ B_{H^*} ~=~ \Pi/\cI.
   \end{equation}
We now make a few observations about $B$.  
   \begin{itemize}\label{page:B}
   \item[(i)] An argument similar to the proof of Lemma \ref{lem:dividing by p} shows that $B$ inherits from $\Pi$ the structure of a pre-$\psi^p$-algebra.
   \item[(ii)] For any $\delta \in \Delta$ and $k \geq 1$, the image $\overline{\psi^{pk}f_\delta}$ of the element $\psi^{pk}f_\delta$ in $\Gr^*B \otimes \bZ/p$ is $0$.  This is because 
   \[
   \psi^{pk}f_\delta ~=~ p\alpha + \alpha^\prime
   \]
for some elements $\alpha, \alpha^\prime \in B$ with $\alpha$ in some positive filtration and $\alpha^\prime$ in a filtration strictly greater than that of $\alpha$.
   \item[(iii)] The unstable $\cA$-algebra generators and relations for $H^*$ were the only choices we made in the process of defining the pre-$\psi^p$-algebra $B$.  Thus, the construction of the proposed lift $B$ is canonical up to a choice of generators and relations for $H^*$.
   \end{itemize}

\subsubsection{$B_{H^*}$ is a lift of $H^*$}

Now we show that the pre-$\psi^p$-algebra $B = B_{H^*}$ constructed above is in fact a lift of $H^*$.

First of all, the associated graded algebra mod $p$ of $\Pi$ is the polynomial algebra
   \[
   \Gr^*\Pi \otimes \bZ/p ~=~ \bZ/p \, \lbrack  \,\lbrace X_{(\theta, i_1, \ldots, i_n)} \rbrace \, \rbrack
   \]
in which $n \geq 0$, $\theta \in \Theta$, and the $i_k$ are in the range $\eqref{eq:range of i}$.  The generator $X_{(\theta, i_1, \ldots, i_n)}$ has degree $2d_\theta + 2(p-1)\sum_{j=1}^n i_j$.  (We omit the bars above these elements for typographical reasons.)  In view of observation (ii) above, it follows that the associated graded algebra mod $p$ of $B$ is the quotient
   \begin{equation}\label{eq:rep of B}
   \Gr^*B \otimes \bZ/p ~=~ \frac{\bZ/p \, \lbrack \, \lbrace X_{(\theta, i_1, \ldots, i_n)} \rbrace \, \rbrack}{\langle f_\delta \, \colon \, \delta \in \Delta \rangle}.
   \end{equation}
Its operations $P^i$ are given on the generators by the equation
   \[
   P^i X_{(\theta, i_1, \ldots, i_n)} ~=~
   \begin{cases}
   X_{(\theta, i_1, \ldots, i_n)}    & \text{ if } i = 0 \\
   X_{(\theta, i_1, \ldots, i_n, i)} & \text{ if } 1 \leq i \leq d_\theta - 1 + (p-1) \sum_{j=1}^n\, i_j \\
   X_{(\theta, i_1, \ldots, i_n)}^p  & \text{ if } i = d_\theta + (p-1) \sum_{j=1}^n\, i_j \\
   0                                 & \text{ otherwise}.
   \end{cases}
   \]

Combining this description of $\Gr^*B \otimes \bZ/p$ with the identification \eqref{eq:rep of H} of $H^*$, we see that the map
   \begin{equation}\label{eq:B to H}
   \rho_{H^*} \colon \Gr^* B \otimes \bZ/p \to H^*
   \end{equation}
defined on the generators by
   \[
   \rho_{H^*} X_{(\theta, i_1, \ldots, i_n)} ~=~ \rP^{i_n} \cdots \rP^{i_1}x_\theta
   \]
is a graded $\bZ/p$-algebra isomorphism that commutes with the operations $P^i$.  Now since the Adem relations and the identity $P^0 = \Id$ hold in $H^*$, it follows that these two properties also hold in $\Gr^*B \otimes \bZ/p$.  In particular, $B$ is, in fact, a $\psi^p$-algebra and is a lift of $H^*$ via the canonical map $\rho_{H^*}$.

Furthermore, the second assertion in Theorem \ref{thm:lift} is a consequence of the constructions of $B_{H^*}$ and $\rho_{H^*}$ (see, in particular, \eqref{eq:psi p X} -- \eqref{eq:B to H}).

The proof of Theorem \ref{thm:lift} is complete.


\section{Kuhn's Realization Conjecture}\label{sec:kuhn}

Our last result is a $K$-theoretic analogue of a conjecture of N.\ Kuhn.  This is our first attempt to tackle the problem of realizing $\psi^p$-algebras as $K$-theory of spaces.  Here we consider the simpler notion of a $\psi^p$-module.

In \cite{kuhn} Kuhn made an interesting conjecture, the Realization Conjecture, about the size of the mod $p$ cohomology of topological spaces: The mod $p$ cohomology of a space should be either finite as a set or infinitely generated as a module over the mod $p$ Steenrod algebra.  To rephrase it, the conjecture predicts that a finitely generated $\cA$-module that is not finitely generated as a $\bZ/p$-module cannot be realized as the mod $p$ cohomology of a space.  Kuhn verified this conjecture in the case when the Bockstein is zero in sufficiently high degrees \cite[Theorem 0.1]{kuhn}.  Using reduction steps in Kuhn's paper \cite{kuhn}, the Realization Conjecture was proved recently by L.\ Schwartz \cite{schwartz}.

One naturally wonders if there are analogous results concerning the size of spaces in other cohomology theories.  Using Atiyah's Theorem \ref{thm:atiyah} and Kuhn's original result, we establish such an analogue for $K$-theory.  To generalize the result of Kuhn and Schwartz, we first introduce a $K$-theoretic notion which corresponds to a module over the Steenrod algebra.

\begin{definition}\label{def:psi p module}
Let $p$ be a prime.  We define a \emph{$\psi^p$-module} to be an ordered pair $((M, \lbrace I_n \rbrace), \psi^p)$ (or simply $(M,\psi^p)$ or even just $M$) consisting of a filtered abelian group $(M, \lbrace I_n \rbrace)$ (see \S \ref{subsubsec:filtered group}) and a distinguished endomorphism $\psi^p$, which is required to satisfy the following condition: If $\alpha \in M$ lies in filtration, say, $2q$, then there exist elements $\alpha_i$ ($i = 0, \ldots , q$) in filtration $2q + 2i(p-1)$ such that
   \begin{equation}\label{eq:atiyah module}
   \psi^p(\alpha) ~=~ \sum_{i=0}^q\, p^{q-i} \alpha_i.
   \end{equation}
\end{definition}

An expression as in \eqref{eq:atiyah module} is referred to as an \emph{Atiyah formula for} $\alpha$.  For example, the $K$-theory $K(X)$ with its Adams operation $\psi^p$ is a $\psi^p$-module, provided that $X$ is a space whose integral cohomology is concentrated in even dimensions and is torsionfree.

Now we can ask what a $K$-theoretic analogue of a \emph{finitely generated} $\cA$-module is.  The $\cA$-linear multiples of an element in an $\cA$-module are the finite sums of iterated Steenrod operations acting on that element.  Since Atiyah's result above tells us that the Steenrod operations on $H^*(X;\bZ/p)$ come from Atiyah formula decomposition (eq.\ \eqref{eq:atiyah formula}) of $\psi^p$, a $K$-theoretic analogue of $\cA$-linear multiples should involve iterated applications of $\psi^p$ on Atiyah formula.  We arrive at the following $K$-theoretic finiteness condition, which corresponds to $H^*(X;\bZ/p)$ being a finitely generated $\cA$-module.

\begin{definition}\label{def:psi p fg}
Let $(M, \psi^p)$ be a $\psi^p$-module and let $m_1, \ldots , m_n$ be elements in $M$.   We say that $M$ is $\psi^p$-\emph{finitely generated by} $m_1, \ldots , m_n$ if the following condition is true:  There exist Atiyah formulas
   \begin{equation}
   \label{eq:psi p fg}
   \begin{split}
   \psi^p m_* &~=~ \sum_{j_1=0}^{q_*} \, p^{q_* - j_1}m_{(*,j_1)} \quad ( 1 \leq * \leq n) \\
   \psi^p m_{(1,j_1)} &~=~ \sum_{j_2=0}^{q_1 + j_1(p-1)} \, p^{q_1 + j_1(p-1) - j_2}m_{(1,j_1,j_2)}  \quad (0 \leq j_1 \leq q_1) \\
   \vdots     & \\
   \psi^p m_{(n,j_1)} &~=~ \sum_{j_2=0}^{q_n + j_1(p-1)} \, p^{q_n + j_1(p-1) - j_2}m_{(n,j_1,j_2)} \quad (0 \leq j_1 \leq q_n) \\
   \vdots & \\
   \end{split}
   \end{equation}
etc., etc.\  such that $M$ is generated as an abelian group by the elements $m_{(i,j_1, \ldots, j_r)}$ $(1 \leq i \leq n, r \geq 0)$.  The $\psi^p$-module $M$ is said to be $\psi^p$-\emph{finitely generated} if there exists a finite set of elements $m_1, \ldots , m_n$ in $M$ with the above property.
\end{definition}

Having a $K$-theoretic analogue of a finitely generated $\cA$-module, we are now ready for the promised $K$-theoretic analogue of Kuhn's Realization Conjecture, or Kuhn-Schwartz Theorem.  As for the original Kuhn-Schwartz Theorem, we think of the following result as an assertion about non-realization (by ``nice" spaces) of $\psi^p$-modules that are $\psi^p$-finitely generated but whose underlying abelian groups are not finitely generated.

\pagebreak
  \begin{prop}
  \label{thm:K-theoretic}
  Let $X$ be a torsionfree space of finite type whose integral cohomology is concentrated in even dimensions.  If there exists a prime $p$ for which the $\psi^p$-module $(K(X), \psi^p)$ is $\psi^p$-finitely generated, then the underlying abelian group of $K(X)$ must be finitely generated.
  \end{prop}

\begin{proof}
We begin with three reductions.

\emph{Reduction step 1}.  To show that $K(X)$ is a finitely generated abelian group, it suffices to show that its associated graded group
   \[
   \Gr^* K(X) ~=~ H^*(X;\bZ)
   \]
is such.  To see this, first note that $K(X)$ with the topology induced by the filtration ideals
   \[
   I^n ~=~ \ker(K(X) \to K(X_{n-1})), 
   \]
where $X_{n-1}$ denotes the $n$th skeleton of $X$, is Hausdorff; that is, the intersection $\cap_n I^n$ is $0$.  Indeed, an element $\alpha$ in $\cap_n I^n$ is represented by a map 
   \[
   \alpha \colon X \to BU
   \]
whose restriction to each skeleton $X_{n-1}$ is nullhomotopic; that is, \ $\alpha$ is a phantom map from $X$ to $BU$.  But since $H^n(X;\bQ)$ and $\pi_{n+1}BU \otimes \bQ$ cannot be simultaneously nonzero for any integer $n$, there can be no non-nullhomotopic phantom maps from $X$ to $BU$ (see \cite{mcgibbon} for more information about phantom maps).  Therefore, $\alpha$ must be $0$ and so 
   \[
   \bigcap_n I^n ~=~ 0; 
   \]
that is, $K(X)$ is Hausdorff.

Now if $\Gr^* K(X) = H^*(X;\bZ)$ is a finitely generated abelian group, then there exists an integer $N > 0$ such that 
   \[
   H^*(X;\bZ) ~=~ H^{< N}(X;\bZ) \qquad \text{and}\qquad H^{\geq N}(X;\bZ) = 0. 
   \]
It follows that $K(X)$ admits the finite filtration
   \begin{equation}
   K(X) ~=~ I^0 \supset I^1 \supset \cdots \supset I^{N-1} \supset I^N ~=~ \cdots ~=~ \cap _jI^j ~=~ 0.
   \end{equation}
Moreover, in this filtration of $K(X)$ both $I^{N-1} = H^{N-1}(X;\bZ)$ and each successive quotient are finitely generated abelian groups.  Thus an easy (reverse) induction argument implies that $K(X)$ itself is a finitely generated abelian group.

\emph{Reduction step 2}.  To show that the associated graded $\Gr^* K(X) = H^*(X;\bZ)$ is a finitely generated abelian group, it suffices to show that the mod $p$ associated graded group
   \[
   \Gr^* K(X) \otimes \bZ/p ~=~ H^*(X;\bZ/p)
   \]
is a finite dimensional $\bZ/p$-vector space.  This is because of the torsionfree hypothesis on $H^*(X;\bZ)$.

\emph{Reduction step 3}.  By Kuhn's theorem discussed above, to show that $H^*(X;\bZ/p)$ is a finite dimensional $\bZ/p$-vector space it suffices to show that it is a finitely generated $\cA$-module.

Now suppose that $K(X)$ is $\psi^p$-finitely generated by the elements $m_1, \ldots , m_n$.  The image of an element in the mod $p$ associated graded is in general given by the same name with a bar above it.  We claim that the elements $\overline{m}_1, \ldots , \overline{m}_n$ generate $H^*(X;\bZ/p)$ as an $\cA$-module.  What this means is that every element in $H^*(X;\bZ/p)$ can be written as a finite sum of elements of the form
   \begin{equation}
   \rP^{j_r} \cdots \rP^{j_1} \overline{m}_i \qquad (\rP^j = \rSq^{2j} \text{ if } p = 2)
   \end{equation}
which is equal to $\overline{m_{(i,j_1,\ldots,j_r)}}$.  This is, of course, implied by the hypothesis that $K(X)$ is $\psi^p$-finitely generated by $m_1, \ldots , m_n$.

This finishes the proof of the proposition.
\end{proof}

\begin{remark}
As in the case of modules over the Steenrod algebra, purely algebraic counterexamples are easily constructed.  For example, let $p$ be an arbitrary prime, and consider the abelian group 
   \[
   A ~=~ \bigoplus_{n = 0}^\infty\, \bZ \langle x^{p^n} \rangle
   \]
with $x^{p^n}$ in filtration $2p^n$ and the endomorphism 
   \[
   \psi^p(x^{p^n}) ~=~ x^{p^{n+1}}.
   \]
Using a slight modification of Lemma \ref{lem:Atiyah}, it is easy to check that $A$ is a $\psi^p$-module and is $\psi^p$-finitely generated by $\lbrace x \rbrace$, and yet it is not finitely generated as an abelian group.  Thus, Proposition \ref{thm:K-theoretic} tells us that many algebraically allowed $\psi^p$-modules cannot be realized as the $K$-theory of spaces.
\end{remark}




\end{document}